\documentclass[11pt,a4paper]{article}
 
\usepackage{amsmath,amsthm,amssymb}

\usepackage{cite}

\pdfoutput=1

\usepackage{graphics,graphicx}
\usepackage{epsfig}
\usepackage{multicol}
\usepackage{color}
\makeatletter
\@addtoreset{equation}{section}
\makeatother

\usepackage{braket}
\usepackage{bm}

\begin{document} 
\sloppy

\begin{center}
\LARGE{{\bf J{\' a}nos Bolyai, Carl Friedrich Gauss, Nikolai Lobachevsky and the New Geometry: Foreword}}
\end{center}

\begin{center}
\large{L{\' a}szl{\' o} Jenkovszky$^{a,}$\footnote{jenk@bitp.kiev.ua}, Matthew J. Lake$^{b,c,d,e,f}$\footnote{matthewjlake@narit.or.th} and Vladimir Soloviev$^{g,}$\footnote{vladimir.soloviev@ihep.ru}}
\end{center}
\begin{center}
\emph{$^{a}$ Bogolyubov Institute for Theoretical Physics, \\ National Academy of Science of Ukraine, 03143 Kiev, Ukraine\\}
\emph{$^{b}$ National Astronomical Research Institute of Thailand,  \\ 260 Moo 4, T. Donkaew,  A. Maerim, Chiang Mai 50180, Thailand\\}
\emph{$^{c}$ Department of Physics and Materials Science, Faculty of Science, Chiang Mai University, \\ 239 Huaykaew Road, T. Suthep, A. Muang, Chiang Mai 50200, Thailand\\}
\emph{$^{d}$ School of Physics, Sun Yat-sen University, Guangzhou 510275, China\\}
\emph{$^{e}$ Department of Physics, Babe\c s-Bolyai University, \\ Mihail Kog\u alniceanu Street 1, 400084 Cluj-Napoca, Romania\\}
\emph{$^{f}$ Office of Research Administration, Chiang Mai University, \\ 239 Huaykaew Rd, T. Suthep, A. Muang, \mbox{Chiang Mai 50200}, Thailand\\}
\emph{$^{g}$ A. A. Logunov Institute for High Energy Physics, \\ NRC Kurchatov Institute, 142281  Protvino, Moscow Region, Russia\\}
\vspace{0.1cm}
\end{center}







\section{Predecessors} \label{Sec.1}

Nearly 2300 years ago, the Greek mathematician Euclid of Alexandria laid down the basis of the geometry now known from the textbooks and used in everyday life. 
It was based on a number of postulates and axioms. Almost all of them were generally accepted either as obvious, or as necessary for logical consistency. 
A notable exception was the fifth postulate: if a straight line crosses two other straight lines, forming internal angles that equate to less than $2\pi$ radians, then these two lines---if continued infinitely---must cross somewhere on the plane where the lines and angles are located.
     
The fifth postulate has been something of an embarrassment ever since Euclidean geometry was first formulated, although its validity was never seriously questioned before the modern era. 
Moreover, it became a dogma---part of European culture, supported by authorities such as Isaac Newton, Leonardo da Vinci, Galileo Galilei, Johann Kepler, Joseph-Louis (Giuseppe Luigi) Lagrange and Emmanuel Kant, among others. 
    
It was only in the first half of the 19th century that three great men---J{\' a}nos Bolyai, Carl Friedrich Gauss and Nikolai Lobachevsky, abbreviated in alphabetical order as BGL throughout the rest of this article---independently, but almost simultaneously, succeeded in generalising Euclid’s vision. 
By making the infamous the fifth postulate his last, Euclid had already alluded to its deficiency, and a great number of Greek, Arab and Renaissance-era mathematicians tried to prove, disprove, or replace this troublesome aspect of \mbox{classical geometry}. 
    
Before BGL, at least two men came close to introducing a ``new geometry'', successfully surpassing the flat Euclidean paradigm. Girolamo Saccheri, an Italian monk, took a step in the right direction in his ``Logica democratica'' by attempting to formulate the fifth postulate in a consistent way, using reductio ab absurdum (evidence to the contrary). 
His ideas were further elaborated in the paper “Euclides ab omni naevo findicas” (Euclid, free of any shadow), which was published in 1697, the year he died.

More than a century later, in the period between 1807 and 1816, Ferdinand Karl Schweikart, a German lawyer employed in Kharkov (now Ukraine), developed a version of non-Euclidean geometry he called ``Astralische Geometrie'' (Astral geometry), and hinted at the cosmic scales at which departures from Euclidean geometry may be noticeable. 
Schweikart, an amateur mathematician, did not use formulae in his work. 
This was implemented by his nephew, Taurinus, who published his “Geometria prima elementa” (first element geometry) in 1826. 
In this work, the ``log-spherical'' geometry, preceding that of Bolyai and Lobachevsky, was introduced to prove Euclid’s fifth postulate. 
Later, the great 19th century mathematician Carl Friedrich Gauss familiarised himself with Schweikart’s work.

\section{Non-Euclidean Geometry in a ``Nutshell''}  \label{Sec.2}

The striking coincidence of independent discoveries by J. Bolyai, C.F. Gauss and N.I. Lobachevsky, after more than two thousand years of stagnation, may seem almost miraculous. In fact, the ground was prepared by their predecessors, as mentioned in the previous section. 
Still, the new concept they introduced was neither well received, nor accepted by their contemporaries, causing significant drama in the lives of (at least two) of its creators. 

Here, we do not discuss the relative importance, or priority, of the discoverers. 
Instead, we briefly introduce the concept of the new geometry, without reference to the individual contributions of BGL. 
In the next Section, we highlight the important milestones in these contributors' biographies, allowing the reader to judge their priorities for themselves. 
In the new geometry, apart from the ``Euclidean'' parallels $AD || CD$, the existence of an infinite number of other straight lines, never crossing $CD$, is assumed to be possible (see Figure \ref{Figs.1}). 
The new ``parallelism angle'' $\Pi(x)$ (which is equivalent to $\pi/2$ in the Euclidean case) is related to the distance $x$ (see Figure \ref{Figs.2}) by the formula $\cot(\Pi(x)/2) = qx$, in which $q$ is a parameter. 
                                                   
\begin{center}
\begin{figure}[h] 
\includegraphics[width=8cm]{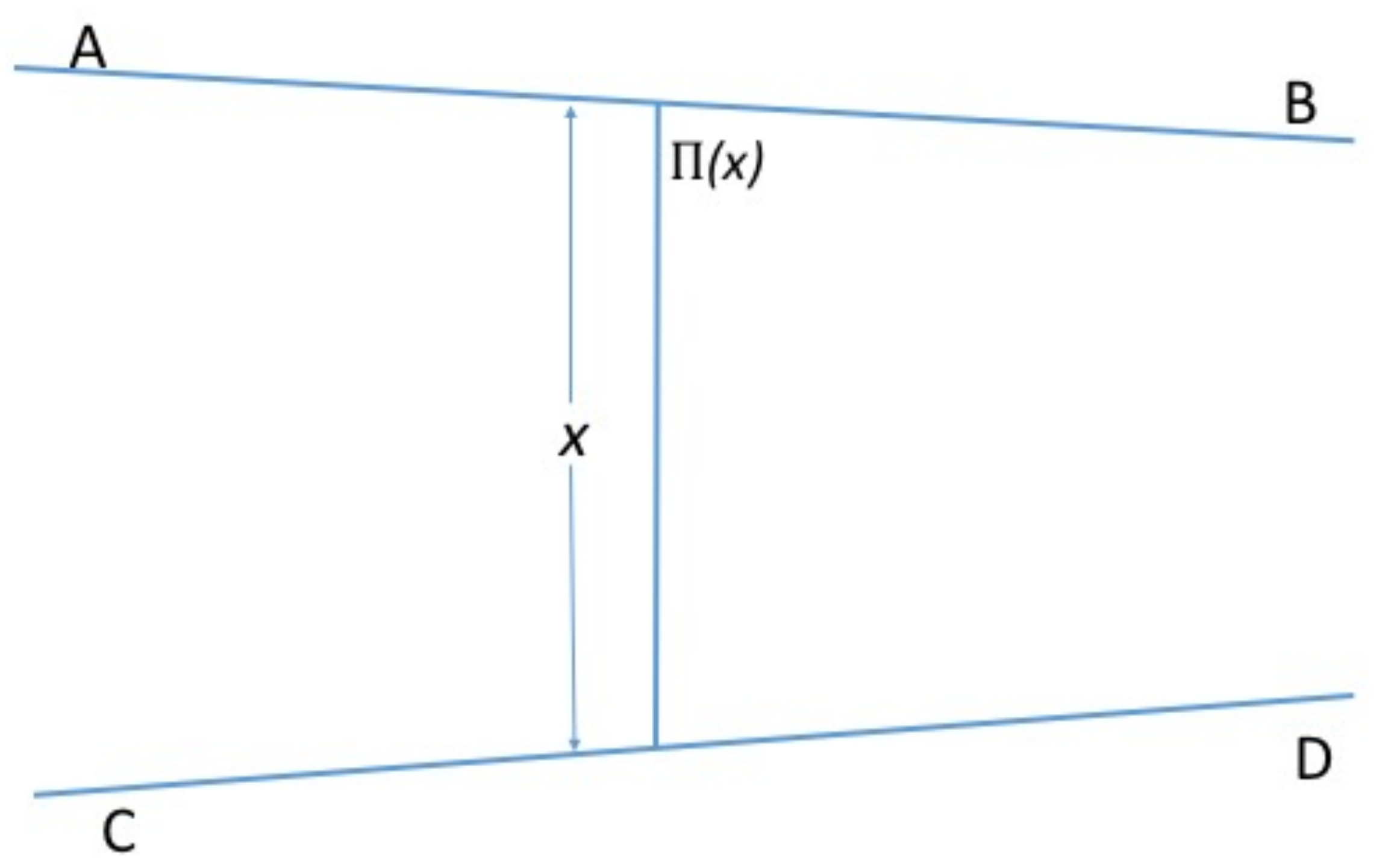}
\caption{In contrast to Euclidean geometry, the sum of the internal angles in the new geometry is no longer $\pi$.\label{Figs.1}}
\end{figure}
\vspace{-11pt}
\end{center}

\begin{center}
\begin{figure}[h] 
\includegraphics[width=8cm]{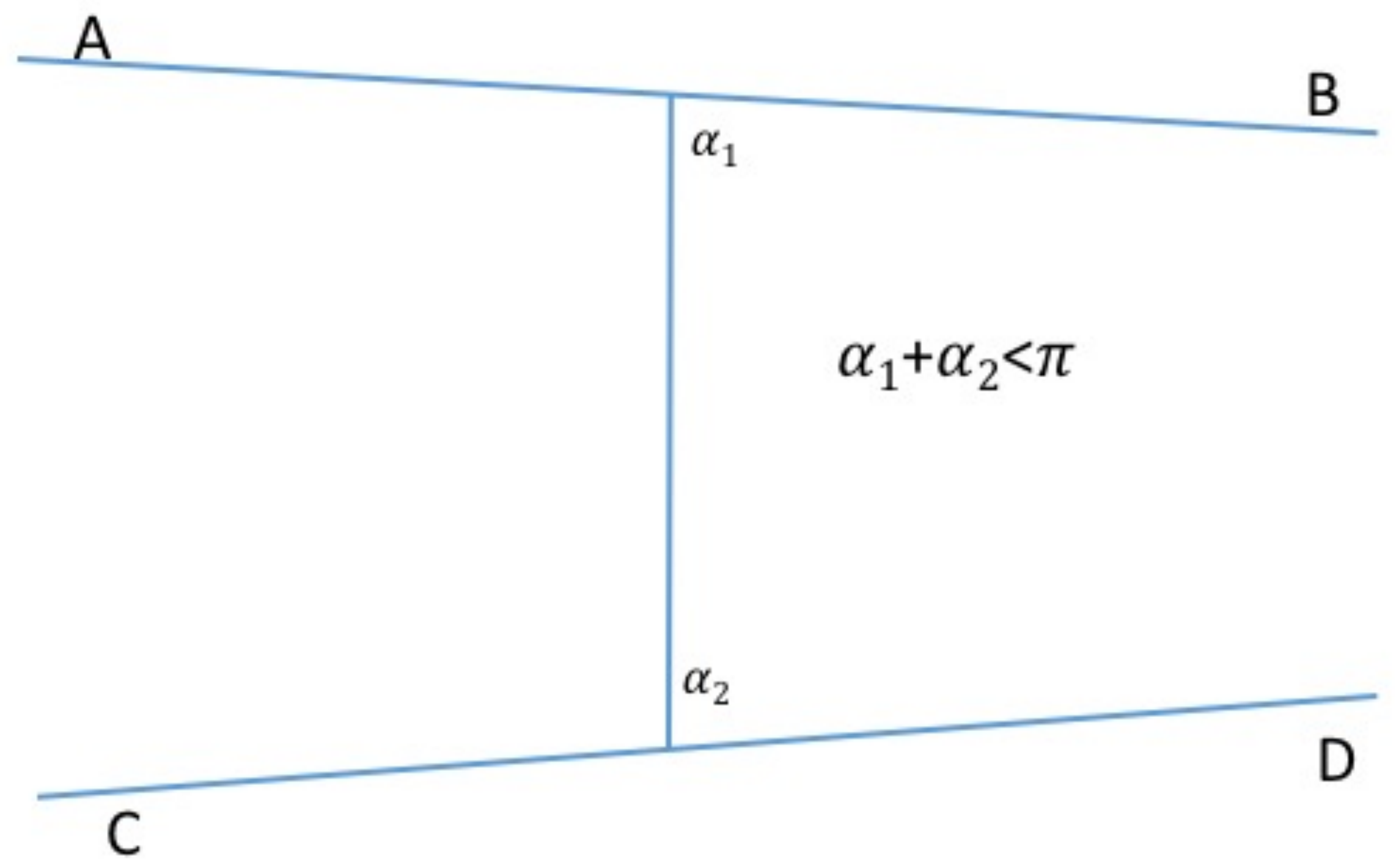}
\caption{The ``parallelism angle'' $\Pi(x)$ is related to the distance $x$.\label{Figs.2}}
\end{figure}
\vspace{-11pt}
\end{center}

The sum of the internal angles in the new geometry is no longer $\pi$: it may be smaller, depending on the length of the sides. 
(See also Figure \ref{Fig.3}, where the curvature on the sides is shown, for pedagogical clarity.)  
J. Bolyai, C.F. Gauss and N.I. Lobachevski were fully aware of the non-observability of any departure from Euclidean geometry within the visible Universe. 
Hence, terms such as ``imaginary'',  ``absolute''  and “pan-“, etc., were at times used to refer to the new geometry. 
The lack of any observable evidence for its existence in nature was, undoubtedly, one of the reasons why their contemporaries were so reluctant to accept it.  
    
\begin{center}
\begin{figure}[h] 
\includegraphics[width=8cm]{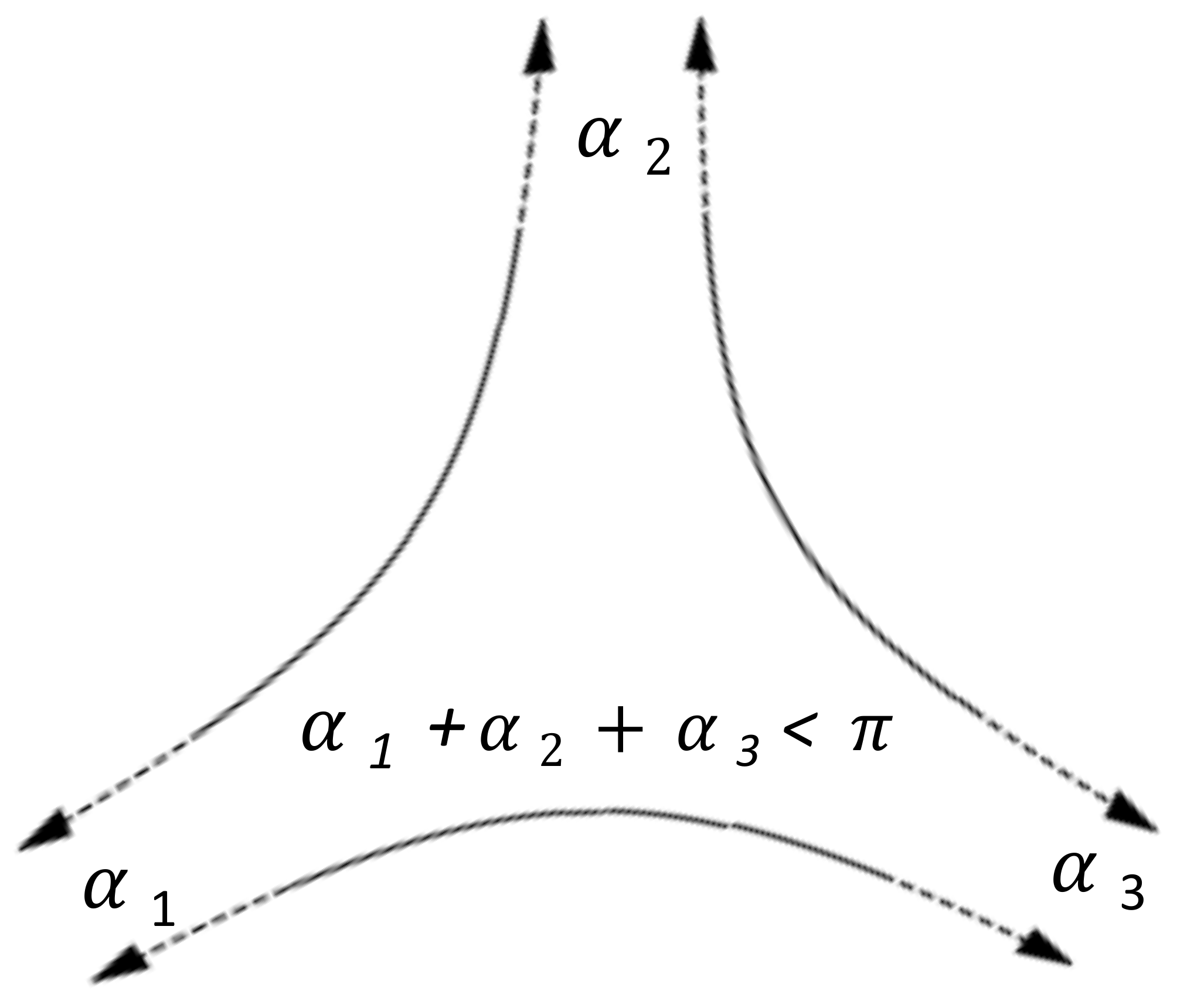}
\caption{Departures from Euclidean geometry may be noticeable only at cosmic distances.\label{Fig.3}}
\end{figure}
\vspace{-11pt}
\end{center}


\section{Carl Friedrich Gauss}  \label{Sec.3}

The eldest member of the BGL trio was born on 30th April 1777, in Braunschweig in the Holy Roman Empire (now Germany). 
In 1795, he entered the University of G{\" o}ttingen, where Farkas Bolyai (J{\' a}nos’ father) became his closest friend. 
After three years of study, Farkas returned to his homeland in Transylvania to teach mathematics in the town of Marosv{\' a}s{\' a}rhely. 
Their friendship survived for decades and generated an extensive exchange of letters, which now provide invaluable information for historians of science.  
    
In 1804, Farkas sent to Gauss his ``proof'' of the fifth postulate. 
In his reply, Gauss indicated an error in his friend’s derivation. 
The work of Gauss was inspired, to a large extent, by a letter he received from Schweikart in 1819, in which the latter explained his ideas. 
Schweikart published his paper on parallel lines, ``Astralische Geometrie'', in 1807, and further developed his ideas from 1812--1816 after moving from Kharkov to Marburg. 
Taurinus continued his uncle’s work in 1824, while in close correspondence with Gauss. 
As mentioned previously, both Schweikart's and Taurino's ideas were close to those of the younger Bolyai. 
However, neither man was able to completely abandon the old Euclidean way of thinking. 
    
After 1816, Gauss, partly inspired by practical interests, began working on geodesics, and in 1828, he published his famous paper on differential geometry (curved surfaces). 
Nonetheless, it was Beltrami who, in 1868, interpreted the new non-Euclidean geometry in terms of surfaces with negative curvature. 
In 1832, Gauss first familiarised himself with the work of J{\' a}nos Bolyai, mediated by Farkas; he later praised Lobachevsky’s “Geometrische Untersuchungen” (Geometric investigations) (see \cite{1}). 
    
Gauss himself did not publish a single paper on non-Euclidean geometry. 
On various occasions---for example, in his private letters---he praised both Lobachevsky and J{\' a}nos Bolyai for their contributions to the development of the new geometry. 
However, he never did so publicly. 
In 1842, Lobachevsky was nominated as a member of the G{\" o}ttingen Scientific Academy, on Gauss’ recommendation. However, his work on geometry \mbox{was ignored}. 

\section{Nikolai Ivanovich Lobachevsky}  \label{Sec.4}

Nikolai Ivanovich Lobachevsky was born in Nizhni Novgorod, on the Volga river, on 20th November 1792. 
However, his studies and career were uniquely connected with the city of Kazan. 
Kazan, another famous city, downstream from Nizhni Novgorod on the Volga, with a new university founded in the first half of the 19th century, was gradually becoming an important regional centre in Eastern Russia. 
The faculty of the new university, invited there from St. Petersburg, and from other important centres of learning in Europe, came with new ideas---in particular, new ideas about geometry, which greatly influenced Lobachevsky in his early years. 
    
The first presentation of his new geometry took place at the Department of Physical and Mathematical Sciences, at Kazan University, on 7th February 1826. 
The subsequent application of the Department, to publish the written version of the presentation entitled “Exposition succinte des principles de la geometrie” (Brief exposition of the principles of geometry), was rejected by the local journal. 
The manuscript was lost. 
A year later, Lobachevsky was elected Rector of Kazan University. 
In 1832, at the age of 40, he married Varvara Musin-Pushkina. 
    
Lobachevsky’s first publication on the new geometry was dated 1829, when the ``Kazanskij vestnik'' published his paper ``On the principles of geometry'' (in Russian). 
In 1832, this paper was submitted for publication to the Russian Academy of Sciences in St. Petersburg. 
The referee’s report, by M.V. Ostrogradsky(!), was totally negative. 
Moreover, in 1834, the magazine ``Syn otechestva'' published a mocking pamphlet, criticising the author and his ideas. 
Lobachevsky shared the fate of many great men and women, whose ideas were ahead of their time. 
    
The overwhelmingly negative reaction of the mathematical community did not shake Lobachevsky’s determination. 
He continued writing and publishing, producing two further notable works: ``G{\' e}om{\' e}trie imaginaire''' (Imaginary geometry, in French) in 1837, and “Geometrische Untersuchungen” (Geometric investigations, in German) in 1840. 
Finally, one year before his death in 1856, while ill and blind, he dictated ``Pangeometria'' \cite{1}, which was published in Russian in 1855, with a French translation following soon after, in 1856.  
    
The Bolyais (both father and son), and Gauss, were familiar with the ``Geometrische Untersuchungen''. 
Farkas, in his 1851 book, and Gauss, in a private letter, both noted and praised it while its author was still alive. 
Nevertheless, during his lifetime, Lobachevsky did not receive the public recognition he deserved. 
He eventually resigned from the position of Rector, partly voluntarily, and partly due to lack of support. 
Following his most active period as a Rector, reformer and lecturer, his career fell into decline. 
Ill and in misery, he died on 12 February  1856. 
     
The serious recognition and revival of his ideas began after his death, but it took another twenty years before his achievements were widely recognised  in Russia, where he had lived and worked for much of his life. 
On the occasion of Lobachevsky's centenary, in 1893, a prize bearing his name was established by the Russian Academy of Science. 
The first five winners were, in chronological order: Sophus Lie (1893), Wilhelm Killing (1900), David Hilbert (1903), Hermann Weil (1917) and Elie Cartan (1937). 

\section{J{\' a}nos Bolyai}  \label{Sec.5}

J{\' a}nos Bolyai was the youngest and, perhaps, the most tragic figure in the dramatic birth of the new geometry. 
Born in Kolozsv{\' a}r (now Cluj-Napoca, Romania) in 1802, he moved with his father Farkas to the town of Marosv{\' a}s{\' a}rhely, also in Transylvania, in 1804. 
His mathematical abilities were recognised from an early age: at the age of 14, he had already mastered differential and integral calculus. 

In 1818, his father asked his friend Gauss to support his son’s studies by inviting him to the University of G{\" o}ttingen, but, for some reason, Gauss was reluctant to do so. 
Consequently, J{\' a}nos entered the Engineering Academy in Vienna. 
After graduation, he was appointed as a custom’s officer at Temesv{\' a}r (now Timisoara). 
In 1820, he informed his father that he had found a way to prove the fifth postulate of Euclid. 
Farkas responded by discouraging his son from spending time on a pursuit that he considered hopeless. 
Farkas played an important role in his son’s life. 
He was, at times, a supporter, adviser, referee and censor to János. 
The Bolyais were descended from an old, highly educated, Hungarian noble family. 
    
Despite his father’s advice, J{\' a}nos continued his search; in 1824, this resulted in the discovery of the mathematical relation between the length of a perpendicular and the angle of its asymptote (see Figures \ref{Figs.1}--\ref{Fig.3}). 
``I have created a world from nothing'', he wrote to his father. 
One can see from his manuscripts, sketches and letters that in 1820, J{\' a}nos was already on the right track by considering the limit of a large circle. 
Sadly, most of Bolyai's works on the new geometry were never published. 
His original manuscripts can be found in various libraries and museums; for example, in Marosv{\' a}s{\' a}rhely.   
    
In February 1825, J{\' a}nos sent a handwritten manuscript on the new geometry to Farkas. 
Unable to understand the new ideas it contained, Farkas searched for deficiencies in his son’s work. 
Finally, in February 1829, he agreed to include his son’s results as an Appendix in his own book, ``Tentamen Juventutem…'', a textbook aimed at introducing mathematics to young people. 
The Appendix (in Latin) was entitled “Appendix scientium spatii absolute veram exhibeus” (\emph{Scientists' absolutely true space, Appendix}).  
The book appeared in 1831. 
A copy was immediately sent, by Farkas, to Gauss, but it did not reach its intended recipient because the area in which Gauss lived was plagued by cholera at the time. 
Another copy reached its intended destination at the beginning of 1832. 
Gauss responded immediately, in March 1832, and his response was disastrous for J{\' a}nos: ``You may be surprised that I will not praise your son’s work since praising it would mean praising myself . . . his ideas almost coincide with my way of thinking already for 30--35 years from now . . . Myself, I also intended to publish these results, but once my friend’s son did this job, I am happy to be free from this obligation'', he wrote. 
     
Gauss' letter to Farkas made the father happy, but not his son. 
J{\' a}nos suspected that his ideas were stolen. 
Furthermore, in 1848, having received Lobachevsky’s ``Geometrische Untersuchungen'', he suspected that ``Lobachevsky'' was a pseudonym, used by Gauss to undermine his priority. 
After this shock, J{\' a}nos started reading Lobachevsky’s work critically, and realised that it was close to what he himself had achieved \cite{2}. 
    
J{\' a}nos continued working hard, extending his research to topology. 
However, deprived of public recognition, lacking a successful professional career, and without critical support from family and friends, he ran into depression. 
Lonely and miserable, he died in Marosv{\' a}s{\' a}rhely in 1860. 
His grave was found by local enthusiasts, with difficulty, in 1893~\cite{3}.
    
Recognition for Bolyai’s achievements also came only after his death. 
Between 1897 and 1904, a series of new editions of the ``Tentamen'' \cite{4}, including the famous ``Appendix'' and marking Bolyai’s centenary, appeared in Budapest. 
Furthermore, the Bolyai Prize was founded by the Hungarian Academy of Sciences. 
Its first winners were Henry Poincar{\' e} (1905) and David Hilbert (1905).  

\section{Successors} \label{Sec.6}

Roughly two decades after the work of Bolyai and Lobachevsky, Georg Friedrich Bernhard Riemann continued investigations of the new geometry. 
In his early years, he planned to become a pastor like his father, studying theology at the University of G{\" ot}tingen. 
Fortunately for modern science, he also studied mathematics under Gauss, who advised him to switch professions. 
Four years later, in 1854, in his profound habilitation lecture ``On the hypotheses which lie at the foundation of geometry'', attended by Gauss, Riemann proposed a far-reaching program for the future development of geometry. 
He intended to generalise Gauss' theory of surfaces to higher dimensions, but he did not succeed in realising this plan. 
There were two likely reasons for this. 
First, he was ahead of his time and risked becoming an outcast by pursuing this program right at the start of his academic career. 
Additionally, despite his brilliant habilitation lecture, he did not succeed in obtaining his desired position of Extraordinary Professor at the university. 
Riemann was able to obtain this position only after four years of successful work in more mainstream areas of contemporary mathematics. 
The text of the talk was not even published during his lifetime. 
Second, and unfortunately for modern science, Riemann's life was too short (1826--1866); he died at 39. 

A large part of his unpublished works were thrown away by his housekeeper. 
Nevertheless, he succeeded in constructing the basis of a new form of mathematics, now called Riemannian geometry. 
The Euclidean and non-Euclidean geometries of BGL may be considered its special cases. 
Among Riemann's many achievements are the first definition of a manifold of arbitrary dimension, the idea of taking the metric tensor as a basic characteristic and defining feature of a given geometry, and the definition of the curvature tensor. 
Manifolds with both positive and negative curvature, which may vary from point to point, were admissible objects in this geometry. 
The Riemann curvature tensor appeared in his 1861 article devoted to thermodynamics, which was submitted to a competition announced by the Parisian Academy. 
In the small part of this work which dealt with geometry, Riemann specified the constraints required for the metric to be reducible to the Euclidean (Pythagorean) form.
        
The next major development occurred when Elwin Bruno Christoffel introduced the notion of covariant differentiation, and the entities were later named after him as Christoffel symbols. 
He also expressed the curvature tensor by means of these symbols. 
For this reason, this tensor is often called Riemann--Christoffel tensor. 
Later, other brilliant mathematicians, including Felix Klein, Gregorio Ricci-Curbastro, Tullio Levi-Civita and Luigi Bianchi, contributed to the development of geometry along similar lines. 
Levi-Civita introduced the concept of the parallel displacement of vectors, along arbitrary curves in a manifold.

In the 20th century, the major revolution in geometry was intimately connected with its applications in physics. 
In 1907, Hermann Minkowski proclaimed that, after the discovery of Special Relativity by Lorentz, Einstein, and Poincar{\' e}, the Newtonian concepts of absolute space and time should be replaced by a new higher-dimensional object---four-dimensional space-time. 
According to Minkowski, ``the world in space and time, in a sense, is a four-dimensional, non-Euclidean manifold'' demonstrating ``the greatest triumph applied mathematics has ever shown''. 
In his memory, this space-time is now called Minkowskian, or pseudo-Euclidean.

Finally, in 1915, a curved generalisation of Minkowski space, now referred to as pseudo-Riemannian geometry, became a crucial ingredient in modern physics, as the geometric basis of the new theory of gravity, ``General Relativity'', created by Albert Einstein. 
Crucially, Einstein's work was supported, and indeed enabled, by the work of mathematicians, including Marcel Grossmann and David Hilbert. 
From the spring of 1915, Einstein also maintained active correspondence with Levi-Civita. 
The formalism of absolute differential calculus, developed by Ricci-Curbastro and Levi-Civita, was used to realise Einstein's ideas regarding the equivalence of inertial and gravitational mass, based on the geometrisation of the gravitational force. 
Space--time became curved, with test bodies moving along geodesics. 
It was Marcel Grossmann who taught Einstein non-Euclidean geometry, and David Hilbert who found the Lagrangian of General Relativity, from which Einstein's gravitational field equations can be derived via a standard variational procedure.

Einstein was happy when, using his new theory, he obtained the observed value of the precession of the perihelion of Mercury. 
He immediately tried to apply the geometric theory of gravity to the whole Universe, considering the matter distribution to be dust-like, homogeneous, and ``at rest'' with respect to space--time. 
In this way, he obtained a static solution in which the space-like slices of the four-dimensional space--time obeyed the geometry of a three-dimensional sphere. 
The radius of this Universe was inversely proportional to the square root of the average mass density of matter. 
To achieve this, Einstein had to introduce a new force---a repulsive term in the field equations generated  by a new constant that he called the ``cosmological constant''. 

The most striking prediction to follow from the equations of General Relativity was discovered by Alexander Friedmann---the double ``n’' was added to his name by Einstein in his reply to Friedman’s first paper---in 1922. 
As a mathematician, Friedmann considered the field equations in more detail than his predecessors, Einstein and De Sitter, and succeeded in finding a set of dynamical equations for the whole Universe, which are now named after him. 
The Friedmann equations had new solutions describing the expansion (or contraction) of the Universe, both with and without a cosmological constant term.  

The dramatic history of the birth, evolution and present state of non-Euclidean geometry inspired physicists and mathematicians from many countries, including Nikolai Alexandrovich Chernikov (in Dubna), L{\' a}szl{\' o} Jenkovszky (Kiev), Elem{\' e}r Kiss (Marosv{\' a}s{\' a}rhely) and Istv{\' a}n Lovas (Budapest), to organise a series of international conferences dedicated to exploring its history and applications. 
The conferences, called BGL, after Bolyai, Gauss and Lobachevsky, were held alternately at different locations connected with its creators (see~\cite{5}). 
The last one was held in Kiev (Kyiv), Ukraine in 2019.

The 2022 conference,  planned in Lviv, Western Ukraine was postponed because of the war. 
The series will undoubtedly continue as soon as peace has returned.





\end{document}